\documentstyle[12pt,amssymb,twoside,pb-diagram,fullpage]{amsart} 
\newtheorem{claim}{Claim}[section]
\newtheorem{conj}[claim]{Conjecture}
\newtheorem{cor}[claim]{Corollary}
\newtheorem{example}[claim]{Example}
\newtheorem{prop}[claim]{Proposition}
\newtheorem{lemma}[claim]{Lemma}

\title{The Gau\ss--Manin connection on the Hodge structures}
\author{M.Rovinsky}
\begin{document}
\maketitle

The Gau\ss--Manin connection is an extra structure on 
the de Rham cohomology of any algebraic variety, $\nabla:
H^{\ast}_{dR/k}\longrightarrow\Omega_k^1\otimes_kH^{\ast}_{dR/k}$ 
(its definition will appear below). If one believes the Hodge conjecture 
then for a given pure Hodge structure $H$ there is at most one connection 
$\nabla$ such that $H$ is a Hodge substructure of a cohomology group of a 
smooth projective complex variety with $\nabla$ induced by the 
Gau\ss--Manin connection. Independently of the Hodge conjecture, there 
are at most countably many connections $\nabla$ on a given pure Hodge 
structure $H$ such that $H$ is a Hodge substructure of a cohomology 
group of a smooth projective complex variety with $\nabla$ induced 
by the Gau\ss--Manin connection (cf. Corollary \ref{count}).

The original motivation for this paper were the properties of 
the forgetful functor $$\left\{\begin{array}{c} 
\mbox{graded-polarizable Hodge structures equipped with} \\
\mbox{a connection respecting the weight filtration and} \\
\mbox{polarizations, and satisfying the Griffiths transversality}
\end{array}\right\}\stackrel{\Phi}{\longrightarrow}\left\{
\mbox{Hodge structures}\right\},$$
meaning the following three questions: 
\begin{itemize} \item For a complex algebraic variety $X$ is the 
Gau\ss--Manin connection determined by the Hodge structure 
$H^{\ast}(X({\Bbb C}))$? 
\item If yes, does there exist a functor $\Psi$ right inverse to 
$\Phi$ such that for each geometric mixed Hodge structure $H$ the 
pair $\Psi(H)$ coincides with $H$ endowed with the Gau\ss--Manin 
connection, and $\Psi(H\otimes H')\cong\Psi(H)\otimes\Psi(H')$? 
\item If the $\nabla$ is determined by the Hodge 
structure in a unique way, how to express $\nabla$ in terms 
of the Hodge structure? (This assumes that there should be a 
certain supply of functions on Shimura varieties classifying 
the Hodge structures.) \end{itemize}

In general these questions are very difficult. In this paper we consider 
some special cases of the problem. In particular, in Proposition 
\ref{h-t-g-m} below we show that such a functor $\Psi$ exists for the 
restriction of $\Phi$ to the subcategory of mixed Tate structures and 
it is unique. The connection constructed there is non-integrable in 
general, so, if any Hodge substructure of a geometric Hodge--Tate 
structure is again geometric, its integrability gives a non-trivial 
necessary condition for a Hodge--Tate structure to be geometric. 

We also compute explicitly the Gau\ss--Manin connection in terms of 
the Hodge structure for some geometric Hodge stuctures of small rank. 

I would like to thank P.Deligne for pointing out non-integrability 
of the connection constructed in Proposition \ref{h-t-g-m}. 
I am grateful to Andrey Levin for the inspiring discussions, 
to the I.H.E.S. for its hospitality and to the 
European Post-Doctoral Institute for its support. 

\section{The absolute Gau\ss--Manin connection: definition}
Let $X_{\bullet}$ be a smooth simplicial quasiprojective variety over 
a field $k$ of characteristic zero. For any $i\ge 0$ reduction modulo 
the ideal generated by $(i+1)$-forms on the base field $k$ gives the 
following short exact sequences of complexes of sheaves of exterior 
powers of absolute K\"ahler differentials on $X_{\bullet}$ 
\begin{equation} \label{main} 0\longrightarrow
\Omega^{i+1}_k\wedge\Omega^{\ge p-1}_{X_{\bullet}}
\longrightarrow\Omega^i_k\wedge\Omega^{\ge p}_{X_{\bullet}}
\stackrel{[~\cdot~]_i}{\longrightarrow}\Omega^i_k\otimes_k
\Omega^{\ge p}_{X_{\bullet}/k}\longrightarrow 0. \end{equation} 
This gives rise to the sequence of homomorphisms with 
composition $\nabla$ called {\it the Gau\ss--Manin connection}:
\begin{equation} \label{gm-def} \begin{diagram} 
\node{\Omega^i_k\otimes_kH^q_{dR/k}(X_{\bullet})} \arrow[2]{e,t}{\nabla} 
\arrow{s,l}{\cong} \node{} 
\node{\Omega^{i+1}_k\otimes_kH^q_{dR/k}(X_{\bullet})} \\
\node{{\Bbb H}^{q+i}(\Omega^i_k\otimes_k\Omega^{\bullet}_{X_{\bullet}/k})} 
\arrow{e,t}{{\rm coboundary}} \node{{\Bbb H}^{q+i+1}(\Omega^{i+1}_k\wedge
\Omega^{\bullet-1}_{X_{\bullet}})} \arrow{e,t}{[~\cdot~]_{i+1\ast}} 
\node{{\Bbb H}^{q+i}(\Omega^{i+1}_k\otimes_k\Omega^{\bullet-1}
_{X_{\bullet}/k})} \arrow{n,r}{\cong} \end{diagram} \end{equation}

Let $\overline{X}_{\bullet}$ be a smooth compactification of $X_{\bullet}$ 
by a divisor with normal crossings $D_{\bullet}$. Then the de Rham 
cohomology of $\overline{X}_{\bullet}$ is identified canonically 
with the hypercohomology of the simplicial logarithmic de Rham complex 
$\Omega^{\bullet}_{\overline{X}_{\bullet}/k}(\log D_{\bullet})$ on 
$\overline{X}_{\bullet}$, so one can define the Hodge filtration on the de 
Rham cohomology groups by $$F^pH^q_{dR/k}(X_{\bullet})={\rm image}\left(
{\Bbb H}^q(\overline{X}_{\bullet},\Omega^{\ge p}_{\overline{X}_{\bullet}/k}
(\log D_{\bullet}))\longrightarrow H^q_{dR/k}(X_{\bullet})\right).$$ 

One has the following short exact sequences of complexes of sheaves of 
exterior powers of absolute K\"ahler differentials on $\overline{X}
_{\bullet}$ which is a logarithmic version of the sequence (\ref{main}): 
\begin{equation} \label{main-prime} \!\!\!\!\!\!\!\!\!\!\!\!\!\!
0\longrightarrow\Omega^{i+1}_k\wedge\Omega^{\ge p-1}_{
\overline{X}_{\bullet}}(\log D_{\bullet})\longrightarrow
\Omega^i_k\wedge\Omega^{\ge p}_{\overline{X}_{\bullet}}
(\log D_{\bullet})\stackrel{[~\cdot~]_i}{\longrightarrow}
\Omega^i_k\otimes_k\Omega^{\ge p}_{\overline{X}_{\bullet}/k}
(\log D_{\bullet})\longrightarrow 0. \end{equation} 
Replacing the complex $\Omega^{\bullet}_{X_{\bullet}}$ on $X_{\bullet}$ by 
the complex $\Omega^{\ge p}_{\overline{X}_{\bullet}}(\log D_{\bullet})$ on 
$\overline{X}_{\bullet}$ in the commutative diagram (\ref{gm-def}), we get 
the Griffiths transversality property: 
\begin{equation} \label{griffiths} \Omega^i_k\otimes_kF^pH^q_{dR/k}
(X_{\bullet})\stackrel{\nabla}{\longrightarrow}
\Omega^{i+1}_k\otimes_kF^{p-1}H^q_{dR/k}(X_{\bullet}). \end{equation}

\section{Basic properties of the Gau\ss--Manin connection} 
Though the most of what follows is presumably valid for arbitrary smooth 
simplicial schemes, we restrict ourselves to the case of smooth proper 
varieties. 

The following properties of the Gau\ss--Manin connection are almost immediate. 
\begin{itemize} \item As $\nabla$ coincides with the first differential 
in the Leray spectral sequence $E^{s,t}_1=\Omega^s_k\otimes_kH^t_{dR/k}
(X_{\bullet})$ converging to $H^{s+t}_{dR/{\Bbb Q}}(X_{\bullet})$ 
(associated to the filtration by powers of the ideal $\Omega^{\ge 1}_k$ 
in $\Omega^{\bullet}_k$), one has $\nabla^2=0$. 
\item Since each morphism $f:X_{\bullet}\longrightarrow Y_{\bullet}$ 
of smooth simplicial schemes induces a morphism of the Leray spectral 
sequences for $X/k$ and $Y/k$, the Gau\ss--Manin connection is 
functorial with respect to pull-backs. 
\item It is easy to see from the K\"unneth decomposition $\Omega^{\bullet}
_{X\times Y}={\rm pr}_X^{\ast}\Omega^{\bullet}_X\otimes_{{\cal O}_{{X
\times Y}}}{\rm pr}_Y^{\ast}\Omega^{\bullet}_Y$ for a pair of smooth 
proper schemes $X$, $Y$, that the Gau\ss--Manin connection on 
$H^{\ast}_{dR/k}(X\times_kY)$ coincides with the tensor product of the 
Gau\ss--Manin connections on $H^{\ast}_{dR/k}(X)$ and $H^{\ast}_{dR/k}(Y)$. 
\item Combining the latter with the functoriality with respect to pull-backs 
applied to the diagonal embedding, we get the Leibniz rule for cup-products. 
\item In fact, one has a stronger version of the Gau\ss--Manin connection. 
Namely, it is defined also on the entire ``conjugate spectral sequence'' 
$\nabla:E^{a,b}_2=H^a(X,{\cal H}^b_{dR/k})\longrightarrow\Omega^1_k\otimes_k
E^{a,b}_2$. It comes from the composition of morphisms of sheaves:
$${\cal H}^b_{dR/k}={\cal H}^b(\Omega^{\bullet}_{X/k})\longrightarrow
{\cal H}^{b+1}(\Omega^1_k\wedge\Omega^{\bullet-1}_X)\longrightarrow\Omega^1_k
\otimes_k{\cal H}^b_{dR/k},$$
where the first map is the coboundary in the cohomological sequence 
for (\ref{main}) with $i=p=0$ and the second is evident. This 
implies that the coniveau filtration is also respected by the 
Gau\ss--Manin connection. In particular, all elements of the 
subgroups $CH^q(X)/CH^q_{alg}(X)\subset E^{q,q}_2$ are horizontal. 
\item The above properties together with Poincar\'e duality give 
the functoriality with respect to Gysin maps. 
\item There is an increasing filtration $W_{\bullet}$ on the complex 
$\Omega^{\bullet}_{\overline{X}_{\bullet}/k}(\log D_{\bullet})$ given 
by limiting the number of logarithmic poles, and inducing a functorial 
increasing weight filtration $W_{\bullet}$ on the cohomology groups 
$H^q_{dR/k}(X_{\bullet})$. The weight filtration is 
obviously respected by the Gau\ss--Manin connection: 
$\nabla(W_wH^q_{dR/k}(X_{\bullet}))\subseteq
\Omega^1_k\otimes_kW_wH^q_{dR/k}(X_{\bullet})$. 
\item One can also state some properties of the gauge-equivalence class 
of the Gau\ss--Manin connection using the properties of variations 
of Hodge structures. \end{itemize}

\begin{prop} If the Hodge conjecture for $H^{2q}$ is true then for a given 
pure effective Hodge structure $H$ of weight $q$ there is at most one 
connection $\nabla$ such that $H$ is a Hodge substructure of the $q$-th 
cohomology group of a smooth projective complex variety with $\nabla$ 
induced by the Gau\ss--Manin connection. \end{prop} 
{\it Proof.} Suppose that $H$ as a pure Hodge structure is isomorphic 
to Hodge substructures of both $H^q(X)$ and of $H^q(Y)$ for some smooth 
projective complex varieties $X$ and $Y$. By Lefschetz hyperplane section 
theorem we may suppose that $\dim X=q$. Then there is a morphism of 
Hodge structures $\alpha:H^q(X)\longrightarrow H^q(Y)$ commuting with 
embeddings $H\hookrightarrow H^q(X)$ and $H\hookrightarrow H^q(Y)$. 
The class of $\alpha$ is an element in $H^{2q}(X\times Y)(q)$ of Hodge 
type $(0,0)$, and thus, is presented by an algebraic cycle $\gamma$. 
By a standard argument $\gamma$ induces a morphism of pairs 
$(H^q(X),\nabla_X)\longrightarrow(H^q(Y),\nabla_Y)$. \hfill $\Box$

\vspace{5mm}

The image $\Omega_H$ of the ${\Bbb Q}$-linear map $H\otimes H^{\vee}
\stackrel{\langle\nabla~\cdot~,~\cdot~\rangle}{-\!\!\!-\!\!\!-\!\!\!%
\longrightarrow}\Omega^1_{{\Bbb C}}$, which is a ${\Bbb Q}$-subspace 
(of dimension $\le({\rm rk}H)^2$) in $\Omega^1_{{\Bbb C}}$, is one of 
basic invariants of the connection $\nabla_H$. Here $H^{\vee}$ is the 
Hodge structure dual to $H$. It follows from the compatibility of 
$\nabla$ with polarizations, duality and tensor products that 
$\Omega_H=\Omega_{H^{\vee}}$, $\Omega_H=\Omega_{H^{\otimes M}}$ 
for any integer $M\ge 1$ and 
$\Omega_{H_1\otimes H_2}\subseteq\Omega_{H_1}+\Omega_{H_2}$. 

If $X$ is a complex algebraic variety and $H=H^q(X)$ is the group of 
singular cohomology of $X({\Bbb C})$ with ${\Bbb Z}$-coefficients for 
some integer $q\ge 0$, one has the de Rham isomorphism $H\otimes{\Bbb C}
\stackrel{\sim}{\longrightarrow}H_{dR/{\Bbb C}}$, and polarizations 
$Q_w:gr^W_wH\times gr^W_wH\longrightarrow{\Bbb Z}(-w)$, 
$Q_w(a,b)=(-1)^wQ_w(b,a)$ for each integer $w$. Denote by $\nabla_H$, or 
simply by $\nabla$, the composition $H\hookrightarrow H\otimes{\Bbb C}
\stackrel{\sim}{\longrightarrow} H_{dR/{\Bbb C}}
\stackrel{\nabla}{\longrightarrow}\Omega^1_{{\Bbb C}}\otimes_{{\Bbb C}}
H_{dR/{\Bbb C}}\stackrel{\sim}{\longrightarrow}\Omega^1_{{\Bbb C}}\otimes H$. 
After a choice of a basis of $H$ one can view this map as a 
$({\rm rk}H\times{\rm rk}H)$-matrix with entries in $\Omega^1_{{\Bbb C}}$. 

The Gau\ss--Manin connection respects the polarizations in the sense 
$$\nabla_{{\Bbb Z}(-w)}Q_w(a,b)=Q_w(\nabla_{gr^W_wH}(a),b)
+Q_w(a,\nabla_{gr^W_wH}(b)).$$

If $w$ is even one can choose such a basis $\{e_1,\dots,e_N\}$ of 
$gr^W_wH$ that the matrix of the polarization is diagonal: $Q_w(e_i,e_j)
=(2\pi\sqrt{-1})^w\lambda_i\delta_{ij}$ for some rational $\lambda_i$'s. 
If $\Omega=(\omega_{ij})$ is the matrix of $\nabla_{gr^W_wH}$ 
in this basis then $\omega_{ii}=-\frac{w}{2}\frac{d\pi}{\pi}$ and
$\omega_{ij}=-\frac{\lambda_i}{\lambda_j}\omega_{ji}$ for $i\neq j$.

If $w$ is odd one can choose a basis $\{e_1,\dots,e_N\}$ of $gr^W_wH\otimes
{\Bbb Q}$ where the matrix of $Q_w$ is equal to $$(2\pi\sqrt{-1})^{-w}
\cdot\left(\begin{array}{cc} 0 & I \\ -I & 0 \end{array}\right)$$
If $\Omega=\left(\Omega_{ij}\right)_{1\le i,j\le 2}$ is the matrix of 
$\nabla_{gr^W_wH}$ in this basis then $\Omega_{12}=\Omega_{12}^t$, 
$\Omega_{21}=\Omega_{21}^t$ and $\Omega_{11}+\Omega_{22}^t=-w
\frac{d\pi}{\pi}\cdot I$. This gives a (very rough) estimate 
$$\dim_{{\Bbb Q}}\Omega_H\le\left\{\begin{array}{ll} \frac{N(N-1)}{2} & 
\mbox{if $H$ is a pure structure of weight 0} \\ \frac{N(N-(-1)^w)}{2}+1 
& \mbox{if $H$ is a pure structure of weight $w$.} \end{array}\right.$$ 

We will need the following 
\begin{prop}[Katz, \cite{nilp-conn}, \cite{hodge-2}] 
\label{hodge-horizontal} Let $S$ be a pathwise connected and locally simply 
connected topological space and $(H,F,W)$ is a family of mixed Hodge 
structures on $S$ such that each of the families $(gr^W_nH,F,W)$ is 
polarizable. Suppose that the Hodge filtration $F$ is locally constant, 
i.e., comes from a filtration of the complexification $H_{{\Bbb C}}$ 
by local subsystems. Then there exists a finite \'etale cover 
$\pi:S'\longrightarrow S$ such that $\pi^{\ast}(H,F,W)$ is a constant 
family of mixed Hodge structures on $S'$. \hfill $\Box$ \end{prop}

\begin{lemma}[{\bf Rigidity}] \label{Rigidity} Let $\overline{f}:
{\frak X}\longrightarrow T$ be a morphism of complex smooth algebraic 
varieties. Denote by $f$ the natural morphism of topological spaces 
${\frak X}({\Bbb C})\longrightarrow T({\Bbb C})$. Suppose that 
$R^qf_{\ast}{\Bbb Z}$ is a local system of isomorphic Hodge structures. 

Then each path $\gamma:[0,1]\longrightarrow T({\Bbb C})$ gives rise to 
an isomorphism $$((R^qf_{\ast}{\Bbb Z})_{\gamma(0)},\nabla_{\gamma(0)})
\stackrel{\sim}{\longrightarrow}((R^qf_{\ast}{\Bbb Z})_{\gamma(1)},
\nabla_{\gamma(1)}).$$ \end{lemma}
{\it Proof.} Since the Hodge filtration on the stalks of 
$R^qf_{\ast}{\Bbb Z}$ defines a filtration of $R^qf_{\ast}{\Bbb C}$ 
by local subsystems, the Proposition \ref{hodge-horizontal} 
implies that the local system $R^qf_{\ast}{\Bbb Z}$ becomes 
trivial on a finite cover of $T$. Without a loss of generality 
we replace $T$ with such a cover, and moreover, assume $T$ 
affine and connected. 

Then for any point $s$ of $T({\Bbb C})$ the restriction $H^q_{dR/{\Bbb C}}
({\frak X})\stackrel{r_s}{\longrightarrow}H_{dR/{\Bbb C}}^q({\frak X}_s)$ 
is surjective and its kernel is independent of $s$. Applying the 
functoriality of the Gau\ss--Manin connection to $r_s$, one gets 
that $\nabla(\ker r_s)\subset\Omega^1_{{\Bbb C}}\otimes\ker r_s$, 
and thus, $\left((R^qf_{\ast}{\Bbb Q})_s,\nabla_s\right)\cong\left(
H^q({\frak X})/\ker r_s,\nabla_{{\frak X}}\right)$. \hfill $\Box$

\begin{cor} \label{count} For a given Hodge structure $H$ there is 
at most a countable set of $\nabla$'s such that $(H,\nabla)$ is a 
cohomology group of a smooth projective complex variety with the 
Gau\ss--Manin connection. \end{cor}
{\it Proof.} There exists such a countable set of families of smooth
projective complex varieties that each variety is isomorphic to an 
element of at least one of the families.\footnote{For each pair of integers 
$1\le r<N$, a collection of integers $2\le d_1\le d_2\le\dots\le d_r$, an 
integer $s\ge 0$ and a collection of elements $P_1,\dots,P_s\in Sym^{\bullet}
(Sym^{d_1}{\Bbb Q}^N\oplus\dots\oplus Sym^{d_r}{\Bbb Q}^N)$ one defines a 
family by homogeneous equations of degrees $d_1,d_2,\dots,d_r$ whose 
coefficients are zeroes of polynomials $P_1,\dots,P_s$.} It follows from 
\cite{som}, that the fibers of the period map are algebraic, and thus, we 
can apply Lemma \ref{Rigidity} to conclude the proof. \hfill $\Box$

\begin{prop} \label{alg-number} Suppose that for a smooth proper 
algebraic variety $X$ over ${\Bbb C}$ and an integer $q\ge 0$ the 
Hodge filtration on $H_{dR/{\Bbb C}}^q(X_)$ is horizontal, i.e., 
$\nabla F^jH_{dR/{\Bbb C}}^q(X)\subseteq\Omega^1_{{\Bbb C}}
\otimes_{{\Bbb C}}F^jH_{dR/{\Bbb C}}^q(X)$. Then $H^q(X)$ 
is the Hodge structure of the $q$-th cohomology group of a 
variety defined over $\overline{{\Bbb Q}}$. \end{prop}
{\it Proof.} Choose a smooth surjective morphism ${\frak X}
\stackrel{\pi}{\longrightarrow}S$ of varieties over $\overline{{\Bbb Q}}$ 
and a generic point $s_0\in S({\Bbb C})$, i.e., an embedding of fields 
$\overline{{\Bbb Q}}(S)\hookrightarrow{\Bbb C}$, 
such that the fiber of ${\frak X}$ over $s_0$ is $X$. 

Then, shrinking $S$, if necessary, we get from Proposition 
\ref{hodge-horizontal}, that the variation $R^q\pi_{\ast}{\Bbb Z}$ of mixed 
Hodge structures is a locally constant local system in \'{e}tale topology. 
Then, by Lemma \ref{Rigidity}, for any point $s_1\in S(\overline{{\Bbb Q}})$ 
we have an isomorphism of mixed Hodge structures equipped with connections 
$\left(H^q({\frak X}_{s_0}({\Bbb C}),\nabla_0\right)\cong\left( 
H^q({\frak X}_{s_1}({\Bbb C}),\nabla_1\right)$. \hfill $\Box$ 

\begin{example} \label{complex-regulators} For any integer $m\ge 2$ those 
elements of $Ext^1_{{\cal HS}}({\Bbb Z},{\Bbb Z}(m))$ corresponding to 
cohomology groups of algebraic varieties, correspond, in fact, to cohomology 
groups of algebraic varieties defined over $\overline{{\Bbb Q}}$. \end{example}
{\it Proof.} In our case $F^{1-m}=\dots=F^0$, so, by the Griffiths 
transversality (\ref{griffiths}), for any integer $1-m\le p\le 0$ one has 
$$\nabla F^p=\nabla F^0\subseteq\Omega^1_{{\Bbb C}}\otimes_{{\Bbb C}}
F^{-1}=\Omega^1_{{\Bbb C}}\otimes_{{\Bbb C}}F^p.$$ 
This means, that the assumptions of the 
Proposition \ref{alg-number} are the case. \hfill $\Box$

\begin{prop} \label{hor-sub} If $H$ is a geometric pure Hodge structure 
then $(H\otimes{\Bbb C})^{\nabla}\otimes_{\overline{{\Bbb Q}}}{\Bbb C}$ 
coincides with $H'\otimes{\Bbb C}$ for a Hodge substructure $H'$ with 
horizontal Hodge filtration. \end{prop}
{\it Proof.} This follows from the exactness of the sequence 
$$H^{\ast}_{dR/{\Bbb Q}}(X)\longrightarrow H^{\ast}_{dR/{\Bbb C}}(X)
\stackrel{\nabla}{\longrightarrow}\Omega^1_{{\Bbb C}}\otimes_k
H^{\ast}_{dR/{\Bbb C}}(X)$$ shown in Proposition 4 of \cite{ep}. \hfill $\Box$ 

\begin{conj} \label{trans} For any geometric Hodge structure $H$ 
the horizontal subspace $H^{\nabla}$ is a Hodge substructure 
isomorphic to a power of ${\Bbb Q}(0)$. \end{conj}

If $H^{\nabla}$ is a Hodge substructure of weight $w$, there is a 
non-degenerate pairing $H^{\nabla}\otimes H^{\nabla}\longrightarrow
{\Bbb Q}(w)$ compatible with the connections, so there are no rational 
horizontal elements in geometric Hodge structures of non-zero weight. 

By Propositions \ref{alg-number} and \ref{hor-sub} a horizontal Hodge 
structure comes from a variety defined over $\overline{{\Bbb Q}}$, 
and thus the Conjecture \ref{trans} is equivalent to transcendence 
of certain periods of varieties over $\overline{{\Bbb Q}}$. 

\section{Hodge--Tate structures} \label{hodge-tate}
\subsection{Calculation of the connection on the logarithmic 
structures}
\label{Cotcotls} Consider the first relative cohomology group of 
${\Bbb G}_m$ modulo $\{1,a\}$ for some $a\in k$. To calculate
this group we present ${\Bbb G}_m$ as the complement of 
${\Bbb P}^1$ to the divisor $(0)+(\infty)$ and then 
$$H_{dR}^1({\Bbb G}_m,\{1,a\})={\Bbb H}^1({\cal O}_{{\Bbb P}^1}(-(1)-(a))
\longrightarrow\Omega^1_{{\Bbb P}^1}((0)+(\infty))).$$ 

We can choose a covering ${\Bbb P}^1=U_0\cup U_1$, say, with
$U_0={\Bbb P}^1\backslash\{a\}$ and 
$U_1={\Bbb P}^1\backslash\{1\}$. 

1-cocycles in the \v{C}ech--de Rham complex are 
collections $(f_{ij},\omega_i)$ with $i,j\in\{0,1\}$ and 
$df_{ij}=\omega_i-\omega_j$, where 
$$f_{ij}\in{\cal O}(U_{01})={\Bbb C}\left[\frac{z-1}{z-a},\frac{z-a}{z-1}
\right]\quad\mbox{and}\quad
\omega_i\in\Omega^1_{{\Bbb P}^1/{\Bbb C}}((0)+(\infty))(U_i).$$ 
Note, however, that adding a coboundary, we may assume the 1-forms 
$\omega_i$ to be regular at the points 1 and $a$, and therefore 
the 1-form $df_{ij}$ to be also regular at the points
1 and $a$. Since $df_{ij}$ is regular everywhere on
the projective line, $df_{ij}=0$, so $f_{ij}$ is a constant
and $\omega_0=\omega_1$. 
The latter means that 
$$\omega_0=\omega_1\in F^1H^1_{dR/{\Bbb C}}({\Bbb G}_m,\{1,a\})=
\Gamma({\Bbb P}^1,\Omega^1_{{\Bbb P}^1/{\Bbb C}}((0)+(\infty)))
=\left\langle\frac{dz}{z}\right\rangle_{{\Bbb C}}.$$

Finally, 
$$H^1_{dR/{\Bbb C}}({\Bbb G}_m,\{1,a\})=
\left\{\left(b,c\frac{dz}{z}\right)|b,c\in{\Bbb C}\right\},$$
where $(1,0)$ denotes the 1-cocycle presented by the function 
1 on $U_{01}$. 

The group $H^1_{dR/{\Bbb C}}({\Bbb G}_m,\{1,a\})$ fits into the exact
sequence $$0\longrightarrow H^0_{dR/{\Bbb C}}({\Bbb G}_m)\longrightarrow 
H^0_{dR/{\Bbb C}}(\{1,a\})\longrightarrow H^1_{dR/{\Bbb C}}({\Bbb G}_m,
\{1,a\})\longrightarrow H^1_{dR/{\Bbb C}}({\Bbb G}_m)\longrightarrow 0,$$
where $H^0_{dR/{\Bbb C}}(\{1,a\})={\Bbb C}\oplus{\Bbb C}$, the first map is 
diagonal and the second map is given by $(s,t)\longmapsto(s-t,0)$. 
In particular, denote by $e_0$ the image of $(1,0)$, equivalently, 
$$e_0=(1,0)\in\check{{\cal C}}^1({\cal O}_{{\Bbb P}^1}(-(1)-(a)))\oplus
\check{{\cal C}}^0(\Omega^1_{{\Bbb P}^1/{\Bbb C}}(-(1)-(a))).$$ 

Note, that $e_0$ lifts tautologically to a 1-cocycle in 
the first term 
$$\check{{\cal C}}^1({\cal O}_{{\Bbb P}^1}(-(1)-(a)))\oplus
\check{{\cal C}}^0(\Omega^1_{{\Bbb P}^1/{\Bbb Q}}
(\log((0)+(a)+(1)+(\infty)))(-(1)-(a)))$$ 
of the absolute \v{C}ech--de Rham complex,
and therefore, $\nabla e_0=0$. 

To calculate $\nabla(\frac{dz}{z})$ we lift the (relative) form 
$\frac{dz}{z}$ to a section $\eta_j$ of the sheaf of absolute 1-forms
$${\cal O}_{{\Bbb P}^1}(-(1)-(a))\otimes_{{\cal O}_{{\Bbb P}^1}}
\Omega^1_{{\Bbb P}^1} \oplus\frac{dz}{z}\cdot{\cal O}_{U_0}\oplus
\frac{d(z/a)}{z/a}\cdot{\cal O}_{U_1}\longrightarrow\!\!\!\!\!\rightarrow
\Omega^1_{{\Bbb P}^1/{\Bbb C}}((0)+(\infty))$$ 
over each element $U_j$ of the covering, say, 
$\eta_0=\frac{dz}{z}$, $\eta_1=\frac{d(z/a)}{z/a}$. 
Then the coboundary of the 1-cochain $(\eta_j)$, a 2-cocycle in 
$$\check{{\cal C}}^1(\Omega^1_{{\Bbb C}}\otimes_{{\Bbb C}}
{\cal O}_{{\Bbb P}^1}(-(1)-(a)))\oplus\check{{\cal C}}^0(\Omega^1_{{\Bbb C}}
\otimes_{{\Bbb C}}\Omega^1_{{\Bbb P}^1/{\Bbb C}}(-(1)-(a))),
\quad\mbox{is}\quad(\frac{da}{a},0).$$ 

This gives the formula $\nabla(\frac{dz}{z})=\frac{da}{a}\otimes e_0$. 

$e_0$ is an integral generator of $W_0H^1({\Bbb G}_m,\{1,a\})$.
Another generator $e_2$ of $H^1({\Bbb G}_m,\{1,a\})$ is
$e_2=[\frac{1}{2\pi i}\frac{dz}{z}]-\frac{\log a}{2\pi i}e_0$ 
for an arbitrary choice of $\log a\in{\Bbb C}$. Then 
$\nabla e_2=\frac{a^{-1}da-d(\log a)}{2\pi i}\!\otimes\!e_0
-\frac{d\pi}{\pi}\!\otimes\!e_2$. 

Finally, for Hodge structure $H$ of rank 2 and weights 0 and 2, 
and any element $\xi$ of $H\cong H^1({\Bbb G}_m,\{1,a\};{\Bbb Q})$ we get 
\begin{equation} \label{log} \nabla\xi=\left(d(\pi z_0)+i
e^{-\pi iz_0}de^{\pi iz_0}\right)\!\otimes\!\frac{e_0}{\pi}
-\frac{d\pi}{\pi}\!\otimes\!\xi, \end{equation} where 
$e_0\neq 0$ is an element of $W_0$ and $\xi-z_0e_0\in F^1$. 
It is easy to see that (\ref{log}) is independent of $e_0$. 

\subsection{The Gau\ss--Manin connection on arbitrary Hodge--Tate structures}
\begin{lemma}\footnote{This is a very particular case of Deligne's 
construction of a functorial splitting of Hodge structures (see, e.g.,
\cite{BZ}, Definition and Proposition 2.6).} \label{h-t-lattice}
For any Hodge--Tate structure $H$ inclusion maps induce the decomposition 
\begin{equation} \label{h-t-decomp} \oplus_kF^k\cap(W_{2k}\otimes{\Bbb C})
\stackrel{\sim}{\longrightarrow}H\otimes{\Bbb C}. \end{equation} \end{lemma}
{\it Proof.} Note that $F^k\cap(W_{2k-2}\otimes{\Bbb C})=0$, so the subspaces 
$F^k\cap(W_{2k}\otimes{\Bbb C})$ and $F^l\cap(W_{2l}\otimes{\Bbb C})$ 
intersect trivially when $k\neq l$, and the canonical projection $\varphi_k:
F^k\cap(W_{2k}\otimes{\Bbb C})\longrightarrow gr^W_{2k}\otimes{\Bbb C}$ is an 
isomorphism. Comparison of dimensions of both sides of (\ref{h-t-decomp})
ensures it is an isomorphism. \hfill $\Box$

As a consequence of this Lemma, we get a new ${\Bbb Q}$-structure 
\begin{equation} \label{h-t-conn-form} 
\oplus_k(2\pi i)^k\varphi_k^{-1}(gr^W_{2k}) \end{equation}
for complexification of arbitrary mixed Tate structure. 

Let ${\frak H}$ be an abelian category of Hodge--Tate structures 
containing all logarithmic structures, invariant under Tate twists 
and containing each Hodge substructure of each its object. 
Consider the category whose objects are objects of ${\frak H}$ equipped 
with a connection satisfying the Griffiths transversality and morphisms 
are morphisms of Hodge structures commuting with connections. 
\begin{prop} \label{h-t-g-m} \begin{itemize} 
\item There is a unique functor $H\longmapsto(H,\nabla_H)$ 
right inverse to the forgetful functor $$\left\{\begin{array}{c} 
\mbox{{\rm objects of} ${\frak H}$ {\rm equipped with a connection}} \\ 
\mbox{{\rm satisfying the Griffiths transversality}} \\ 
\end{array}\right\}\longrightarrow{\frak H}$$ 
such that for a logarithmic structure $H$ the connection $\nabla_H$ 
coincides with the Gau\ss--Manin connection calculated above, 
and $\nabla_H$ induces the same connection on 
$H\otimes{\Bbb C}=H(1)\otimes{\Bbb C}$ as $\nabla_{H(1)}$. 
\item \label{integrability-gen} The above connection on a Hodge--Tate 
structure $H$ is integrable if and only if the above connections 
on $W_{2k}/W_{2k-6}(k-3)$ are integrable for all integer $k$. 
\end{itemize} \end{prop} 
{\it Proof.} It follows from the functoriality, applied to the morphism 
$W_k\longrightarrow H$, that $\nabla$ respects the weight filtration, 
and therefore, from the calculation for the logarithmic structures, that 
the connection on the structure ${\Bbb Z}(0)$ is zero. This implies that 
$\nabla:W_{2k}(k)\longrightarrow\Omega^1_{{\Bbb C}}\otimes W_{2k-2}$. 
Combining these with the Griffiths transversality 
$\nabla:F^k\longrightarrow\Omega^1_{{\Bbb C}}\otimes F^{k-1}$, we get 
$$(2\pi i)^k\varphi^{-1}_k(gr^W_{2k})\stackrel{\nabla}{\longrightarrow}
\Omega^1_{{\Bbb C}}\otimes_{{\Bbb C}}(F^{k-1}\cap(W_{2k-2})_{{\Bbb C}}).$$ 

It is clear from the decomposition (\ref{h-t-decomp}) that it 
suffices to construct the latter maps, and that $\nabla$ is 
integrable if and only if the restrictions of $\nabla^2$ to 
$F^k\bigcap(W_{2k})_{{\Bbb C}}$ vanish for all integer $k$. As 
$F^{k-2}\bigcap(W_{2k-6})_{{\Bbb C}}=0$, the vanishing of the latter 
restrictions is equivalent to the vanishing of the induced maps 
$F^k\bigcap(W_{2k})_{{\Bbb C}}\longrightarrow
\Omega^2_{{\Bbb C}}\otimes(W_{2k}/W_{2k-6})$. 

By the functoriality and compatibility with the Tate twists, 
to construct that map for some $k$ we may identify the space 
$F^k\cap(W_{2k})_{{\Bbb C}}$ with the space $F^1H'\cap(W_2H')_{{\Bbb C}}$, 
and $F^{k-1}\cap(W_{2k-2})_{{\Bbb C}}$ with $F^0H'\cap(W_0H')_{{\Bbb C}}$, 
where $H'=(W_{2k}/W_{2k-4})(k-1)$. In fact, we may suppose that there is 
an exact sequence $$0\longrightarrow{\Bbb Z}(0)^s\longrightarrow 
H'\longrightarrow{\Bbb Z}(-1)\longrightarrow 0$$
for some non-negative integer $s$. Then $H'$ can be identified 
with a Hodge substructure of a direct sum of logarithmic Hodge 
structures, where we have fixed the connection. \hfill $\Box$ 

{\sc Remarks.} 1. It is easy to see that the connection constructed 
in Proposition \ref{h-t-g-m} on the tensor product of two 
Hodge--Tate structures coincides with the tensor product of the 
connections on that Hodge--Tate structures.

2.\label{integrability-3} It follows from Section \ref{Cotcotls} 
that for a Hodge--Tate structure $H$ of rank 3 with weights 0, 2, 4 
the connection $\nabla_H$ is integrable if and only if $e^{\pi iz_0}$ 
and $e^{\pi iz_2}$ are algebraically dependent, where $z_0$ and $z_2$ 
are determined by the conditions $e_0\in W_0\backslash\{0\}$, 
$e_2\in W_2\backslash W_0$, $e_4\in W_4\backslash W_2$, 
$e_2-z_0\cdot e_0\in F^1$ and 
$e_4-z_2\cdot e_2\in F^2(W_4/W_0)_{{\Bbb C}}$. This implies that 
for each element $a\in{\Bbb C}^{{\times}}$ there is a natural 
embedding ${\Bbb C}^{{\times}}/(\overline{{\Bbb Q}(a)})^{{\times}}
\hookrightarrow{\rm Ext}^2({\Bbb Z}(0),{\Bbb Z}(2))$, where 
${\rm Ext}^2$ is calculated in the category of flat Hodge--Tate 
structures. 

\section{Examples: some Hodge structures of rank $\le 3$} Up to the Tate 
twists and the duality the options for a Hodge structure of rank $\le 3$ are: 
\begin{itemize} 
\item it is Hodge--Tate; \item it is pure of weight $0$, or $-1$; 
\item it is an extension of ${\Bbb Q}(0)$ by a pure structure of negative 
weight. \end{itemize}

We have calculated the ``Gau\ss--Manin'' connection for Hodge--Tate 
structures in the previous section, so we eliminate them in what follows. 

\subsection{Pure Hodge structures of rank 2} \label{phsr-2}
For any square-free positive integer $D$ denote by $E_D$ an 
elliptic curve with complex multiplication in ${\Bbb Q}(\sqrt{-D})$. 
Since $\wedge^2H^1(E_D)={\Bbb Q}(-1)$, one has a surjection 
${\cal E}nd H^1(E_D)\stackrel{\sim}{\longrightarrow}H^1(E_D)\otimes H^1(E_D)(1)
\longrightarrow\!\!\!\!\!\rightarrow{\rm Sym}^2H^1(E_D)(1)$
with the kernel ${\Bbb Q}\cdot id$. The classes of endomorphisms of 
the curve generate a Hodge substructure in ${\cal E}nd H^1(E_D)$ 
isomorphic to ${\Bbb Q}(0)\oplus{\Bbb Q}(0)$, so its image in 
${\rm Sym}^2H^1(E_D)(1)$ is isomorphic to ${\Bbb Q}(0)$. 
Let $\alpha\in{\rm Sym}^2H^1(E_D)(1)$ be the element corresponding to 
$1\in{\Bbb Q}(0)$. 

Denote by $M_D^k$ the cokernel of the injection ${\rm Sym}^{k-2}H^1(E_D)(-1)
\stackrel{\cdot\alpha}{\longrightarrow}{\rm Sym}^kH^1(E_D)$. 
This is a geometric Hodge structure of weight $k$, rank 2 
and with Hodge numbers $h^{0,k}=h^{k,0}=1$. 

{\it Alternatively}, one can define $M_D^k(k)$ as a Hodge substructure 
in $H_1(E_D)^{\otimes k}$ as follows. Let $\tau$ be a complex 
multiplication of $E_D$, $\tau_{\ast}$ the induced endomorphism of the 
first homology group, and $\gamma_1$ a non-trivial element of $H_1(E_D)$. 
Set $\gamma_2=\tau_{\ast}\gamma_1$, $\gamma_1^{(k)}=(\gamma_2- 
\overline{\tau}\gamma_1)^{\otimes k}+(\gamma_2-\tau\gamma_1)^{\otimes k}
$ and $\gamma_2^{(k)}=(\tau-\overline{\tau})^{-1}\left((\gamma_2- 
\overline{\tau}\gamma_1)^{\otimes k}-(\gamma_2-\tau\gamma_1)^{\otimes k} 
\right)$. Then $\gamma_1^{(k)}$ and $\gamma_2^{(k)}$ are integral elements, 
generating a Hodge substructure $M^k_D(k)$ with induced polarization 
$$\langle\gamma_1^{(k)},\gamma_1^{(k)}\rangle=(1+(-1)^k)c^k,\quad
\langle\gamma_2^{(k)},\gamma_2^{(k)}\rangle=
-\frac{1+(-1)^k}{(\tau-\overline{\tau})^2}c^k,\quad
\langle\gamma_1^{(k)},\gamma_2^{(k)}\rangle=
\frac{1-(-1)^k}{\overline{\tau}-\tau}c^k,$$
where $c=(\tau-\overline{\tau})\langle\gamma_1,\gamma_2\rangle$.
In particular, when $k$ is even the polarization is of type $x^2+Dy^2$.

\begin{claim} If $\{\gamma_1^{(k)},\gamma_2^{(k)}\}$ is a basis of $M^k_D(k)$ 
with $\gamma_1^{(k)}+(\tau-\overline{\tau})\gamma_2^{(k)}\in F^0$ then 
$$\begin{array}{c} 
\frac{\nabla\gamma_1^{(k)}}{k}=\frac{d\pi}{2\pi}\otimes\gamma_1^{(k)}+
(\tau-\overline{\tau})\left(\frac{d\pi}{\pi}-\frac{w}{4h}
\sum_{n\in({\Bbb Z}/m)^{\times}}\left(\frac{-D}{n}\right)
\frac{d\Gamma(n/m)}{\Gamma(n/m)}\right)\otimes\gamma_2^{(k)}; \\
\frac{\nabla\gamma_2^{(k)}}{k}=\frac{1}{\tau-\overline{\tau}}\left(
\frac{d\pi}{\pi}-\frac{w}{4h}\sum_{n\in({\Bbb Z}/m)^{\times}}
\left(\frac{-D}{n}\right)\frac{d\Gamma(n/m)}{\Gamma(n/m)}\right)
\otimes\gamma_1^{(k)}+\frac{d\pi}{2\pi}\otimes\gamma_2^{(k)}, \end{array}$$
where $w=4$ if $D=1$, $w=6$ if $D=3$, and $w=2$ otherwise; 
$h$ is the number of classes of ideals in the ring of integers in 
${\Bbb Q}(\sqrt{-D})$; and 
$$m=\left\{\begin{array}{ll} D, & \mbox{if $D\equiv 1$ mod 4} \\
4D, & \mbox{otherwise.} \end{array}\right.$$ \end{claim}
{\it Proof.} It will follow from Example \ref{elliptic-example} that $2\pi i
\nabla\gamma_1=(\omega_2d\eta_1-\eta_2d\omega_1)\otimes\gamma_1+(\eta_1
d\omega_1-\omega_1d\eta_1)\otimes\gamma_2$ and $2\pi i\nabla\gamma_2=(\omega_2
d\eta_2-\eta_2d\omega_2)\otimes\gamma_1+(\eta_1d\omega_2-\omega_1d\eta_2)
\otimes\gamma_2$. Rewriting the latter using $\gamma_2=\tau_{\ast}\gamma_1$ as 
$2\pi i\nabla\gamma_2=2\pi i(id\otimes\tau_{\ast})\nabla\gamma_1=\tau
\overline{\tau}(\eta_1d\omega_1-\omega_1d\eta_1)\otimes\gamma_1+(\dots)
\otimes\gamma_2$, we 
get $\eta_2d\omega_1-\omega_2d\eta_1=\overline{\tau}(\eta_1d\omega_1-
\omega_1d\eta_1)$, or equivalently, $d(\eta_2-\overline{\tau}\eta_1)=
(\eta_2-\overline{\tau}\eta_1)\frac{d\omega_1}{\omega_1}$. From the 
Legendre identity $\omega_1\eta_2-\omega_2\eta_1=2\pi i$ and 
$\omega_2=\tau\omega_1$ we get 
$d(\eta_2-\tau\eta_1)=d\left(\frac{2\pi i}{\omega_1}\right)$. 
As a corollary of these formulas, 
$2\pi i\nabla(\gamma_2-\tau\gamma_1)=\omega_1d\left(\frac{2\pi i}{\omega_1}
\right)\otimes(\gamma_2-\tau\gamma_1)$; and $2\pi i\nabla(\gamma_2-
\overline{\tau}\gamma_1)=2\pi i\frac{d\omega_1}{\omega_1}
\otimes(\gamma_2-\overline{\tau}\gamma_1)$.

Finally, $$\nabla\left(\begin{array}{c} \gamma^{(k)}_1 \\
\gamma^{(k)}_2 \end{array}\right)=k\left(\begin{array}{cc} 
\frac{d\pi}{2\pi} & 
\left(\frac{d\omega_1}{\omega_1}-\frac{d\pi}{2\pi}\right)(\tau-
\overline{\tau}) \\
\frac{1}{\tau-\overline{\tau}}\left(\frac{d\omega_1}{\omega_1}-
\frac{d\pi}{2\pi}\right) & \frac{d\pi}{2\pi} \end{array}\right)\left(
\begin{array}{c} \gamma^{(k)}_1 \\ \gamma^{(k)}_2 \end{array}\right)$$

Due to the Chowla--Selberg formula (cf., e.g., \cite{weil}) $\omega_1$ 
is an algebraic multiple of \\
$\pi^{3/2}\prod_{n=1}^{m-1}\Gamma(n/m)^{-w\chi(n)/4h}$,
where $w$ is the number of roots of unity in ${\Bbb Q}(\sqrt{-D})$, $h$ is 
the number of classes of ideals in the ring of integers in 
${\Bbb Q}(\sqrt{-D})$, $m$ is the discriminant of 
${\Bbb Q}(\sqrt{-D})$ and $\chi(n)=\left(\frac{-m}{n}\right)$ is the Jacobi 
symbol. \hfill $\Box$

\begin{itemize} \item Now, a pure polarized Hodge structure of rank 2 of 
even weight $w$ is a triplet of a lattice $H$ of rank 2 with a symmetric 
definite ${\Bbb Q}(-w)$-valued bilinear form on it, and an integer $k>w/2$. 
Then the Hodge filtration is given by
$H_{{\Bbb C}}=F^{w-k}\supset F^k=F^{w-k+1}\supset F^{k+1}=0$, where $F^k$ 
is the isotropic line, on which the corresponding hermitian form is 
positive, so for a fixed triplet of integers $(k>w/2,D>0)$ there is at 
most one isogeny class of polarized Hodge structures of rank 2, weight 
$w$ with polarization of discriminant $D$. On the other hand, the Hodge 
structure $M^{2k}_D(k-w/2)$ gives the example. 
\item To determine a polarized Hodge structure of rank $2$ and odd weight 
$w$ means to fix a lattice $H$ of rank $2$ with an isomorphism 
$\wedge^2H\stackrel{\sim}{\longrightarrow}{\Bbb Z}(-w)$, an integer $k>w/2$ 
and a line $F^k\subset H_{{\Bbb C}}$, on which the corresponding hermitian 
form is positive. Then the Hodge filtration is given 
by $H_{{\Bbb C}}=F^{w-k}\supset F^k=F^{w-k+1}\supset F^{k+1}=0$. 

If $k>\frac{w+1}{2}$ then the Griffiths transversality condition and
Proposition \ref{alg-number} imply that such geometric Hodge
structures correspond to varieties defined over number fields.
At least some of the examples of them are constructed in 
\cite{mod-forms}. In particular, the Hodge structure 
$M^{2k-1}_D(k-\frac{w+1}{2})$ gives the example. 
\item Tensoring with ${\Bbb Z}(\frac{w+1}{2})$, we reduce the case 
$k=\frac{w+1}{2}$ to the case of a Hodge structure $H$ of rank 2 with Hodge 
numbers $h^{0,-1}=h^{-1,0}=1$. Then we have an embedding $H\hookrightarrow 
H_{{\Bbb C}}/F^0$. Denote the one-dimensional ${\Bbb C}$-space 
$H_{{\Bbb C}}/F^0$ by $L$ and by $L^s$ its $s$-th tensor power if $s>0$ and 
the dual of $L^{-s}$ otherwise. Let $\zeta=\zeta_H:L-H\longrightarrow L^{-1}$ 
be the meromorphic function given by 
$$\zeta_H(z)=\frac{1}{z}+\sum_{\lambda\in H\backslash\{0\}}\left(
\frac{1}{z-\lambda}+\frac{1}{\lambda}+\frac{z}{\lambda^2}\right)=z^{-1}-
\sum_{n\ge 2}z^{2n-1}\left(\sum_{\lambda\in H\backslash\{0\}}\lambda
^{-2n}\right).$$ $$\mbox{Set}\quad A=-15\sum_{\lambda\in H\backslash\{0\}}
\lambda^{-4}\in L^{-4}\quad\mbox{and}
\quad B=-35\sum_{\lambda\in H\backslash\{0\}}\lambda^{-6}\in L^{-6}.$$

The Hodge structure $H$ is naturally isomorphic to $H_1(E({\Bbb C}))$ for the 
elliptic curve $E({\Bbb C})=L/H$. Our nearest aim is to obtain an algebraic 
equation (\ref{weier}) of $E$, and therefore, express a basis (\ref{basis}) of 
the de Rham cohomology of $E$ in terms of $H$. This enables us to calculate 
the connection in (\ref{con-dr}) and in Example \ref{elliptic-example}. Set 
$$\begin{array}{c} X=\wp_H(z)=\frac{1}{z^2}+\sum_{\lambda\in 
H\backslash\{0\}}\left(\frac{1}{(z-\lambda)^2}-\frac{1}{\lambda^2}\right)
=z^{-2}-\frac{1}{5}Az^2-\frac{1}{7}Bz^4+O(z^6)\in L^{-2}, \\
Y=-\frac{1}{2}\wp'_H(z)=\frac{1}{2}\zeta''_H(z)=\frac{1}{z^3}+
\sum_{\lambda\in H\backslash\{0\}}\frac{1}{(z-\lambda)^3}
=z^{-3}+\frac{1}{5}Az+\frac{2}{7}Bz^3+O(z^5)\in L^{-3}. \end{array}$$

Since any holomorphic elliptic function is constant, 
we land on the well-known relation
\begin{equation} \label{weier} Y^2=X^3+AX+B. \end{equation}

Fix a twelve-th root $\Delta\in L^{-1}$ of $-4A^3-27B^2$ and define the 
complex numbers $a=A/\Delta^4$ and $b=B/\Delta^6$, so $4a^3+27b^2=-1$. 
Set $x=X/\Delta^2$ and $y=Y/\Delta^3$. 

Also, any connection $\nabla$ on $L^{-2}$ induces conections 
on $L^{-4}$ and $L^{-6}$. We set \begin{equation} \varkappa=
\frac{da}{18b}=-\frac{db}{4a^2}=\frac{2A\nabla B-3B\nabla A}{
2\Delta^{10}}\in\Omega^1_{{\Bbb C}/{\Bbb Q}}. \end{equation}

Starting with elementary identities in $\Omega^1_{{\Bbb C}(x,y)/{\Bbb C}}$ 
\begin{equation} \label{exa} \begin{array}{c} -d\left(\frac{1}{y}\right)=
\frac{(3x^2+a)dx}{2y^3};\quad d\left(\frac{x}{y}\right)
=\frac{(2b+ax-x^3)dx}{2y^3};\quad d\left(\frac{x^2}{y}\right)
=\frac{(x^4+3ax^2+4bx)dx}{2y^3},\end{array} \end{equation} 
and fixing the following differentials of the second kind on $E$ 
\begin{equation} \label{basis} \omega=\frac{dx}{2y}=-dz\quad\quad
\mbox{{\rm and}}\quad\quad\varphi=\frac{xdx}{2y}=d\zeta(z),\end{equation} 
we get then the following congruences modulo exact forms 
$$\omega=\frac{(x^3+ax+b)dx}{2y^3}\equiv\frac{(2ax+3b)dx}{2y^3}; \qquad
\varphi=\frac{(x^4+ax^2+bx)dx}{2y^3}\equiv\frac{(2a^2-9bx)dx}{6y^3}.$$
This information is enough to find the Gau\ss--Manin connection:
$$\nabla\omega=-\left[\frac{dy^2\wedge dx}{4y^3}\right]=
-\left[\frac{(xda+db)\wedge dx}{4y^3}\right]=
\left[\varkappa\wedge\frac{4a^2\cdot dx-18b\cdot xdx}{4y^3}\right].$$
So we get 
\begin{equation} \label{con-dr} \nabla\omega=3\varkappa\otimes[\varphi];
\quad\mbox{{\rm and similarly,}}\quad
\nabla\varphi=a\varkappa\otimes[\omega].\end{equation}

Then from the identities 
$$\langle\nabla\omega,\gamma\rangle+\langle\omega,\nabla\gamma\rangle=d\langle
\omega,\gamma\rangle\qquad\mbox{{\rm and}}\quad\langle\nabla\varphi,\gamma
\rangle+\langle\varphi,\nabla\gamma\rangle=d\langle\varphi,\gamma\rangle,$$ 
$$\mbox{{\rm we get}}\quad\langle\omega,\nabla\gamma\rangle=d\langle\omega,
\gamma\rangle-3\langle\varphi,\gamma\rangle\varkappa\quad\mbox{{\rm and}}
\quad\langle\varphi,\nabla\gamma\rangle=d\langle\varphi,\gamma\rangle-
a\langle\omega,\gamma\rangle\varkappa.$$

Fix such a basis $\{\gamma_1,\gamma_2\}$ of the lattice $H\subset L$ that the 
imaginary part of $\gamma_2/\gamma_1$ is positive. Let 
$\omega_1=\Delta\cdot\gamma_1$ and $\omega_2=\Delta\cdot\gamma_2$ be the 
complex numbers corresponding to $\gamma_1$ and $\gamma_2$, respectively, and 
$\eta_j=\Delta^{-1}\cdot 2\zeta_H(\gamma_j/2)$. 

Using the Legendre relation $\omega_2\eta_1-\omega_1\eta_2=2\pi i$, 
one easily verifies that the following holds. 
\begin{example} \label{elliptic-example} Fix a basis $\{\gamma_1,\gamma_2\}$ 
of $H$ as above. Then, in the above notations, 
\begin{multline*} 2\pi i\nabla\gamma_1=\left[\omega_2d\eta_1-\eta_2d\omega_1
+(a\omega_1\omega_2-3\eta_1\eta_2)\varkappa\right]\otimes\gamma_1 \\ 
+\left[\eta_1d\omega_1-\omega_1d\eta_1+(3\eta_1^2-a\omega_1^2)\varkappa
\right]\otimes\gamma_2; \end{multline*} 
\vspace{-8mm}
\begin{multline*} 2\pi i\nabla\gamma_2=\left[\omega_2d\eta_2-
\eta_2d\omega_2+(a\omega_2^2-3\eta_2^2)\varkappa\right]\otimes\gamma_1 \\ 
+\left[\eta_1d\omega_2-\omega_1d\eta_2+(3\eta_1\eta_2-a\omega_1\omega_2)
\varkappa\right]\otimes\gamma_2. \end{multline*} \end{example} \end{itemize}

\subsection{Pure Hodge structures of rank 3}
Since any pure Hodge structure of odd rank (3 in our case) is of even 
weight, we may suppose the weight is 0. To determine such a polarized 
Hodge structure one has to fix a lattice $H$ of rank 3 with a 
non-degenerate symmetric ${\Bbb Z}$-valued bilinear form on $H$, a 
positive integer $k$ and an isotropic line $F^k\subset H_{{\Bbb C}}$, 
on which the associated hermitian form is positive. 
Then the Hodge filtration is given by 
$H_{{\Bbb C}}=F^{-k}\supset F^0=F^{1-k}\supset F^k=F^1\supset F^{k+1}=0$, 
where $F^0$ is the two-dimensional subspace orthogonal to $F^k$. 
In particular, for a fixed $k$ the classifying space of such structures 
is uniformized by a connected component of the complement of a smooth 
conic in the complex projective plane to the set of its real points. 

\begin{itemize} \item When $k=1$ the Hodge structures ${\frak sl}(H^1(E))$ 
for arbitrary elliptic curves $E$ give all desired Hodge structures with 
polarization of type $x^2-y^2-z^2$, and thus, in this case the calculation 
of the connection is again reduced to the case of elliptic curves. 
\item If $k\ge 2$ then the Griffiths transversality condition and 
Proposition \ref{alg-number} imply that such geometric Hodge structures 
correspond to varieties defined over number fields. \end{itemize}

\subsection{Hodge structures of rank 3 and weights $-k$ for $k\ge 2$ and $0$}
Let for some $t>k/2$ the Hodge filtration on weight-$(-k)$ part be in 
the following range: $F^{-t}(W_{-k})_{{\Bbb C}}=(W_{-k})_{{\Bbb C}}$ and 
$0=F^{t-k+1}(W_{-k})_{{\Bbb C}}\subset F^{t-k}(W_{-k})_{{\Bbb C}}
=F^{1-t}(W_{-k})_{{\Bbb C}}\neq 0$.

\begin{itemize} \item

If $t>k+1$, or $k/2+1\le t<k-1$, then the Hodge filtration 
is obviously horizontal. 

If $t=k+1$ then the filtration looks like 
$0=F^2\subset F^1\subset F^0=F^{-k}\subset F^{-k-1}=H_{{\Bbb C}}$. 
Obviously, $\nabla F^0\subset\Omega^1_{{\Bbb C}}\otimes_{{\Bbb C}}F^0$.
On the other hand, $F^1\subset(W_{-k})_{{\Bbb C}}$, and therefore, 
one has $\nabla F^1\subset\Omega^1_{{\Bbb C}}\otimes_{{\Bbb C}}(F^0\cap
(W_{-k})_{{\Bbb C}})=\Omega^1_{{\Bbb C}}\otimes_{{\Bbb C}}(F^1\cap
(W_{-k})_{{\Bbb C}})=\Omega^1_{{\Bbb C}}\otimes_{{\Bbb C}}F^1$, so
Hodge filtration is again horizontal. 

If $t=k$ then the filtration 
looks like $0=F^1\subset F^0=F^{1-k}\subset F^{-k}=H_{{\Bbb C}}$. 
Obviously, $\nabla F^0\subset\Omega^1_{{\Bbb C}}\otimes_{{\Bbb C}}F^0$.

Finally, we get that for $t\ge k$ and $k/2+1\le t<k-1$ the 
Hodge structures correspond to varieties over $\overline{{\Bbb Q}}$. 
\item If $t=\frac{k+1}{2}$ for an odd $k\ge 5$, then the connection
is trivial on the preimage of ${\Bbb Z}$ under the map 
$F^0\longrightarrow gr^W_0\otimes{\Bbb C}={\Bbb C}$, while the rest is a 
Hodge structure of rank 2, where the connection is already calculated. 
\item If $t=k-1$ for $k\ge 3$ then suppose for simplicity that 
$W_{-k}$ is a Hodge substructure of maximal width of 
$H^{k-2}(E^{k-2})(k-1)$ for an elliptic curve $E$. 
Then it is expected that classes of geometric extensions of 
${\Bbb Q}(0)$ by $H^{k-2}(E^{k-2})(k-1)$ are in one-to-one correspondence 
with certain classes of motivic cohomology: 
$$Ext^1_{{\cal M}{\cal M}}({\Bbb Q}(0),H^{k-2}(E^{k-2})(k-1))=
Gr^1_FCH^{k-1}(E^{k-2};k-1)_{{\Bbb Q}}.$$
As the images of regulators on $K_{\ge 2}$ are countable 
(cf. Claim 2.3.4 of \cite{regul}), such extensions form a countable set. 
\end{itemize}

\subsection{Hodge structures of rank 3 and weights $-1$ and $0$}
These ones are extensions of the Tate structure of weight 0 
by Hodge structures of rank 2 and weight $-1$. 
\begin{itemize}
\item

Let for some $t\ge 2$ the Hodge filtration be in the following range: 
$0=F^t\subset F^{t-1}=F^1\subset F^0=F^{1-t}\subset F^{-t}=H_{{\Bbb C}}$. 

If $t\ge 3$ then the Hodge filtration is horizontal for trivial reasons. 

If $t=2$ then $0=F^2\subset F^1\subset F^0=F^{-1}\subset F^{-2}=H_{{\Bbb C}}$.
In particular, $\nabla F^0\subset\Omega^1_{{\Bbb C}}\otimes_{{\Bbb C}}F^0$.
On the other hand, $F^1\subset(W_{-1})_{{\Bbb C}}$, and therefore,
$$\nabla F^1\subset\Omega^1_{{\Bbb C}}\otimes_{{\Bbb C}}
(F^0\cap(W_{-1})_{{\Bbb C}})=\Omega^1_{{\Bbb C}}\otimes_{{\Bbb C}}
(F^1\cap(W_{-1})_{{\Bbb C}})=\Omega^1_{{\Bbb C}}\otimes_{{\Bbb C}}F^1,$$
so the Hodge filtration is again horizontal. 

This means that for any $t\ge 2$ these Hodge structures correspond
to varieties over $\overline{{\Bbb Q}}$.
\item Now suppose that $t=1$. We have calculated the Gau\ss--Manin connection 
on the weight-$(-1)$ part, the rest is as follows. We keep the notations of 
\S\ref{phsr-2} above with the only exception: we denote by $W_{-1}$ what was 
$H$ there. Fix an element $\gamma\in H$ projecting to the generator $1$ of 
the Tate structure ${\Bbb Z}(0)$. Subtracting from it an element in $F^0$ 
projecting to the same element of the Tate structure ${\Bbb Z}(0)$, we get a 
well-defined element $Z_0$ of the space $L$. Clearly, there is a one-to-one 
correspondence between elements of the torus $L/W_{-1}$ and extensions of the 
Tate structure ${\Bbb Z}(0)$ by the Hodge structure $W_{-1}$. 

Consider the rational 1-form $\frac{y+\beta}{x-\alpha}\cdot\frac{dx}{y}=
\frac{\wp'(z)-2\beta}{\wp(z)-\alpha}dz=\left(-\frac{2}{z}+O(z)\right)dz$ 
from the subspace $F^1H^1_{dR/{\Bbb C}}(E\backslash\{0,Q\})$ with 
logarithmic poles at $Q=(\alpha,\beta)$ and $\infty$, where 
$\alpha=\wp(z_0)$ and $\beta=-\frac{1}{2}\wp'(z_0)$. 

$$\nabla\left(\frac{y+\beta}{x-\alpha}\frac{dx}{y}\right)=\nabla\left(
\left(1+\frac{\beta}{y}\right)\frac{d(x-\alpha)}{x-\alpha}\right)
=\frac{1}{\beta}\left[\left(\frac{y^2d\beta^2}{2y^3}-
\frac{\beta^2dy^2}{2y^3}\right)\wedge
\frac{d(x-\alpha)}{x-\alpha}\right].$$

One checks directly, that $$y^2d\beta^2-\beta^2dy^2\equiv(x-\alpha)
((x^2+\alpha x+\alpha^2+a)d\beta^2-\beta^2d(3\alpha x+a))$$
modulo the multiples of $d(x-\alpha)$, \footnote{Namely, it becomes 
transparent from the identities $$y^2d\beta^2-\beta^2dy^2=(y^2-\beta^2)
d\beta^2-\beta^2d(y^2-\beta^2)=(x-\alpha)(x^2+\alpha x+\alpha^2+a)d\beta^2
-\beta^2d((x-\alpha)^3+(x-\alpha)(3\alpha x+a)).$$} and therefore,
$$\nabla\left(\frac{y+\beta}{x-\alpha}\frac{dx}{y}\right)=\frac{1}{\beta}
\left[\frac{(x^2+\alpha x+\alpha^2+a)d\beta^2-\beta^2d(3\alpha x+a)}{2y^3}
\wedge d(x-\alpha)\right].$$
After some more work\footnote{With a help of the first of the congruences 
(\ref{exa}) we rewrite $\nabla\left(\frac{y+\beta}{x-\alpha}\frac{dx}{y}
\right)$ further as $$\begin{array}{c} 
\left[\frac{(\alpha x+\alpha^2+\frac{2}{3}a)d\beta^2-\beta^2d(3\alpha x+a)}
{2\beta y^3}\wedge d(x-\alpha)\right] \\
=\left[\frac{((\alpha x+\alpha^2+\frac{2}{3}a)(3\alpha^2+a)-
3\beta^2(x+\alpha))d\alpha+(\alpha^2x-\frac{a\alpha}{3}-b)da
+(\alpha x+\alpha^2+\frac{2a}{3})db}{2\beta y^3}\wedge dx\right] \\
\!=\!\frac{(\frac{2}{3}a^2-3\alpha\cdot b-x
(2a\alpha+3b))d\alpha+(\alpha^2x-\frac{a\alpha}{3}-b)da+
(\alpha x+\alpha^2+\frac{2}{3}a)db}{2\beta y^3}\!\wedge\! dx\! \\
=\!\frac{((\frac{2a^2}{3}-3bx)-\alpha
(2ax+3b))d\alpha+(18b(\alpha^2x-\frac{a\alpha}{3}-b)-4a^2
(\alpha x+\alpha^2+\frac{2a}{3}))\varkappa}{2\beta y^3}\!\wedge\! dx \\
=\frac{1}{\beta}\left[d\alpha\wedge\varphi-\alpha d\alpha\wedge\omega
+\frac{(2a+3\alpha^2)(6bx-\frac{4}{3}a^2)-(a\alpha+3b)(4ax+6b)}{2y^3}
\varkappa\wedge dx\right]. \end{array}$$}
we get 
$$\nabla\left(\frac{y+\beta}{x-\alpha}\frac{dx}{y}\right)
=\left[\frac{d\alpha}{\beta}\wedge\varphi-\frac{\alpha d\alpha}{\beta}
\wedge\omega-\frac{2}{\beta}(a\alpha+3b)\varkappa\wedge\omega-
\frac{2}{\beta}(2a+3\alpha^2)\varkappa\wedge\varphi\right].$$

Using an evident expansion $\zeta(z)-\zeta(z-z_0)=\frac{1}{z}+\zeta(z_0)+O(z)$ 
of the elliptic function in $z$ for small $z$, and notation 
$\sigma=\sigma_{W_{-1}}(z)=z\prod_{\lambda\in H\backslash\{0\}}
(1-\frac{z}{\lambda})e^{z/\lambda+z^2/2\lambda^2}$, we get 
$$\frac{y+\beta}{x-\alpha}\cdot\frac{dx}{y}=
2(\zeta(z-z_0)-\zeta(z)+\zeta(z_0))dz=2\cdot d\log\left(
\frac{\sigma(z-z_0)}{\sigma(z)}e^{\zeta(z_0)\cdot z}\right).$$

Obviously, one has $t\sigma_{W_{-1}}(z)=\sigma_{W_{-1}}(tz)$ for any complex 
automorphism $t$ of $W_{-1}$. In particular, it is an odd function in $z$. 
One easily checks that 
$$\sigma(z+\mu)=-e^{(2z+\mu)\zeta(\mu/2)}\sigma(z)\quad
\mbox{for any $\mu\in W_{-1}\backslash 2W_{-1}$}.$$
Namely, second derivatives of logarithm of both sides coincide as they 
are values of the elliptic function $-\wp$ at points $z+\mu$ and $z$,
respectively. This implies that difference of first derivatives of
logarithm of two sides is a constant. Evaluation at $z=-\mu/2$ shows 
that first derivatives of logarithm of both sides coincide. 
This implies that quotient of two sides is a constant. 
Evaluation at $z=-\mu/2$ shows that the equality holds. 

As a corollary, we see that the class $\gamma\in 
H^1(E\backslash\{0,(\alpha,\beta)\}; {\Bbb Z}(1))$ is presented by the 1-form 
$$\frac{y+\beta}{\alpha-x}\frac{dx}{2y}-\zeta(z_0)\cdot\omega-
z_0\cdot\varphi,\quad\mbox{so finally, we get the following}.$$
\begin{example} For any generator $\gamma$ as before, the connection is 
determined by the conditions $\nabla\gamma\in\Omega^1_{{\Bbb C}}\otimes 
W_{-1}$ and 
\begin{multline*} \langle\nabla\gamma,\gamma_j\rangle=\omega_jd\zeta(z_0)-
\eta_j\frac{d\alpha}{2\beta}-\omega_j\frac{\alpha d\alpha}{2\beta}-\eta_jdz_0 
\\ 
+\left(az_0\omega_j-3\zeta(z_0)\eta_j+\frac{2a+3\alpha^2}{\beta}\eta_j-
\frac{a\alpha+3b}{\beta}\omega_j\right)\varkappa.~~\Box\end{multline*} 
\end{example} \end{itemize}

\vspace{5mm}

Current address:
$$\begin{array}{l} \mbox{Independent University of Moscow} \\
\mbox{121002 Moscow} \\ \mbox{B.Vlasievsky Per. 11} 
\\ \mbox{{\tt marat@@mccme.ru}} \end{array}
\quad\!\!\!\!\!\!\!\!\!\!\!\!\!\!\!\!\!\mbox{and}\quad
\begin{array}{l}\mbox{Institute for Information Transmission Problems} \\ 
\mbox{Russian Academy of Sciences} \\ \mbox{Moscow}\end{array}$$
\end{document}